\documentclass[12pt,a4paper]{article}
\usepackage[english]{babel}
\usepackage{amsmath}
\usepackage{amssymb}
\usepackage{amsfonts}
\usepackage{hyperref}

\begin{document}


\renewcommand{\refname}{References}
\renewcommand{\contentsname}{Contents}

\begin{center}
{\huge Global solvability of 1D equations\\ of viscous compressible multi-fluids}
\end{center}

\medskip

\begin{center}
{\large Alexander Mamontov,\quad Dmitriy Prokudin\footnote{The research was supported by the Ministry of Education and Science of the Russian Federation (grant 14.Z50.31.0037).}}
\end{center}

\medskip

\begin{center}
{\large August 25, 2017}
\end{center}

\medskip

\begin{center}
{
Lavrentyev Institute of Hydrodynamics, \\ Siberian Branch of the Russian Academy of Sciences\\
pr. Lavrent'eva 15, Novosibirsk 630090, Russia}
\end{center}

\medskip

\begin{center}
{\bfseries Abstract}
\end{center}


\begin{center}
\begin{minipage}{110mm}
We consider the model of viscous compressible multi-fluids with multiple velocities. We review different formulations of the model and the existence results for boundary value problems. We analyze crucial mathematical difficulties which arise during the proof of the global existence theorems in 1D case.
\end{minipage}
\end{center}

\bigskip


\newpage

\tableofcontents

\bigskip

\section{Introduction}

The description of the motion of multi-component media is an interesting and rather little-studied problem both in physics/mechanics and in mathematics. There is no standard approach to simulating these motions, nor is there any developed mathematical theory concerning the existence, uniqueness and properties of solutions of initial-boundary value problems arising in this simulation.

In the present paper, we choose one of the numerous versions of simulating the motion of multi-component fluid mixtures, namely, a homogeneous mixture of viscous compressible fluids and a multi-velocity model. This means that all components (constituents) of the mixture are present at any point of the space and with the same phase, and each of them has its own local velocity. The interaction between the components occurs via the viscous friction and the exchange between momenta, and also using the heat exchange (in heat-conducting models). This kind of mixtures is also called a multi-fluid, see \cite{mamprok1d.france} for the details.

\section{Mathematical model of multi-fluids}

The mathematical model of a multi-fluid which consists of $N\geqslant 2$ components, includes (see \cite{mamprok1d.Nigm} and \cite{mamprok1d.Raj}) the continuity equations for each constituent
\begin{equation}\label{mamprok1d.21061710}\partial_{t}\rho_{i}+{\rm div}(\rho_{i}\boldsymbol{u}_{i})=0,\quad i=1,\ldots,N,\end{equation}
the momentum equations for each constituent
\begin{equation}\label{mamprok1d.21061711}\partial_{t}(\rho_{i}\boldsymbol{u}_{i})+{\rm div}(\rho_{i}\boldsymbol{u}_{i}\otimes\boldsymbol{u}_{i}) +\nabla p_{i}-{\rm div}\,
{\mathbb S}_{i}=\rho_{i}\boldsymbol{f}_{i} +\boldsymbol{J}_{i},\quad i=1,\ldots,N,\end{equation}
and the energy equations. Here $\rho_i$ is the density of the $i$-th constituent of the multi-fluid, $\boldsymbol{u}_{i}$ is the velocity field, $p_i$ is the pressure,
${\mathbb S}_{i}$ is the viscous part of the stress tensor ${\mathbb{P}}_{i}=-p_i\mathbb{I}+\mathbb{S}_{i}$:
$${\mathbb S}_{i}=\sum\limits_{j=1}^{N}\left(\lambda_{ij}({\rm div}\,\boldsymbol{u}_{j}){\mathbb I}+2\mu_{ij}{\mathbb D}(\boldsymbol{u}_{j})\right),\quad i=1,\ldots,N,$$
where $\lambda_{ij}$ and $\mu_{ij}$ are the viscosity coefficients, ${\mathbb D}$ is the rate of deformation tensor, and ${\mathbb I}$ is the identity tensor. The viscosity coefficients $\lambda_{ij}$ and $\mu_{ij}$ compose the matrices $\boldsymbol{\Lambda}=\{\lambda_{ij}\}_{i, j =1}^{N}$ and $\textbf{M}=\{\mu_{ij}\}_{i, j =1}^{N}$. Finally,\linebreak  $\boldsymbol{f}_{i}=(f_{i1}, \ldots, f_{in})$ denotes the external body force ($n$ is the dimension of the flow), and
$$\boldsymbol{J}_{i}=\sum\limits_{j=1}^{N}a_{ij}(\boldsymbol{u}_{j}-\boldsymbol{u}_{i}),\quad i=1,\ldots,N,\quad a_{ij}=a_{ji},\quad i, j=1,\ldots,N$$
stands for the momentum exchange between the constituents (the momentum supply).

The viscosity matrices in real multi-fluids take non-trivial forms and, generally speaking (see \cite{mamprok1d.mamprok17} for the details), they depend on the concentrations $\displaystyle \xi_{i}=\frac{\rho_{i}}{\rho}$, where $\displaystyle \rho=\sum\limits_{i=1}^{N}\rho_{i}$ is the total density of the multi-fluid. The dependence of the viscosity matrices on the concentrations constitutes serious difficulties. We accept the simplifying assumption of constant viscosities, but with the preservation of necessary properties of positiveness. The thing is that it is very important (physically and mathematically) to validate Second law of thermodynamics which means
\begin{equation}\label{mamprok1d.2106171}\sum\limits_{i=1}^{N}{\mathbb S}_{i}:{\mathbb D}(\boldsymbol{u}_{i})\geqslant 0.\end{equation}
Besides, in order to provide important mathematical property of ellipticity, it is necessary to validate the condition
\begin{equation}\label{mamprok1d.2106172}\sum\limits_{i=1}^{N}\int\limits_{\Omega}{\mathbb S}_{i}:{\mathbb D}(\boldsymbol{u}_{i})d\boldsymbol{x}\geqslant C\sum\limits_{i=1}^{N}
  \int\limits_{\Omega} |\nabla\otimes \boldsymbol{u}_{i}|^2 d\boldsymbol{x},\end{equation}
where $\Omega$ is the flow domain, and $\boldsymbol{u}_{i}|_{\partial\Omega}=0$,  $i=1,\ldots,N$. The formulated positiveness or coercivity can be provided by the following
properties of viscosity matrices: the properties $n\boldsymbol{\Lambda}+2\textbf{M}\geqslant 0$, $\textbf{M}\geqslant 0$ provide \eqref{mamprok1d.2106171}, and the properties
$\boldsymbol{\Lambda}+2\textbf{M}>0$, $\textbf{M}>0$ provide \eqref{mamprok1d.2106172}.

A very important observation is that the viscosity matrices should not be diagonal. Momentum supply $\boldsymbol{J}_{i}$ in the momentum equations gives lower order terms (physically important, but mathematically causing no difficulties), and if the matrices are diagonal then $\boldsymbol{J}_{i}$ is the only connection between the constituents, so we have $N$ Navier--Stokes systems connected only via lower order terms. Earlier such problems were relevant (even in 1D) (see, e.~g., \cite{mamprok1d.kazhpetr78}, \cite{mamprok1d.petr82} and \cite{mamprok1d.zlotn95}), but nowadays such results almost automatically come from the theory of mono-fluid systems governed by the compressible Navier--Stokes equations.
If viscosity matrices are ``complete'' (general) then we have interesting mathematical problems (see the reviews in \cite{mamprok1d.france}, \cite{mamprok1d.semr17} and the citations there).

The crucial problem in the multi-D flows is that an automatic extension of the theory of compressible Navier--Stokes equations to the theory of multi-fluids requires
$${\rm div}\,{\rm div}\,\mathbb S_{i}={\rm const}_{i}\cdot\Delta{\rm div}\,\boldsymbol{u}_{i},\quad i=1,\ldots,N,$$
but we have
$$\left( \begin{array}{c}{\rm div}\,{\rm div}\,\mathbb S_{1} \\ \ldots \\ {\rm div}\,{\rm div}\,\mathbb S_{N} \end{array} \right)=\textbf{N}
\left( \begin{array}{c}\Delta{\rm div}\,\boldsymbol{u}_{1} \\ \ldots \\ \Delta{\rm div}\,\boldsymbol{u}_{N} \end{array} \right),$$
where $\textbf{N}=\{\nu_{ij}\}_{i, j =1}^{N}=\{\lambda_{ij}+2\mu_{ij}\}_{i, j =1}^{N}$ is the matrix of total viscosities. It is possible to obtain results in the case of triangular matrices $\textbf{N}$ (see, e.~g., \cite{mamprok1d.smj12} and \cite{mamprok1d.izvran}). However, for a general~$\textbf{N}$, it is important to radically review the existing techniques developed for the compressible Navier--Stokes system. If $n=1$ then interesting problems also appear (see below).

The authors found a way to overcome the crucial problem of total viscosity matrix structure, that is to introduce the assumptions that
\begin{itemize}
  \item the pressures in all constituents are equal to each other,
  \item the velocities of the constituents in the material derivative operators are replaced by the average velocity.
\end{itemize}
Both assumptions are physically realistic in many situations (see \cite{mamprok1d.mamprok17} for the details), and the main mathematical effect is that the model does not lose richness inherent to the multi-fluid models (different densities and velocities of the constituent are preserved), moreover, the variety of the model grows because this trick allows to eliminate the restrictions of the viscosity matrix, and hence to admit all possible viscous terms.

As a result, we come (see \cite{mamprok1d.mamprok17} and \cite{mamprok1d.semr17}) to the following system of equations:
\begin{equation}\label{mamprok1d.continit}\partial_{t}\rho_{i}+{\rm div\,}(\rho_{i}\boldsymbol{v})=0,\quad i=1, \ldots, N,\end{equation}
\begin{equation}\label{mamprok1d.mominit}\partial_{t}(\rho_{i}\boldsymbol{u}_{i})+{\rm div\,}(\rho_{i}\boldsymbol{v}\otimes\boldsymbol{u}_{i})+\alpha_i\nabla p
={\rm div\,}{\mathbb S}_{i}+\rho_{i}\boldsymbol{f}_{i},\quad i=1, \ldots, N.\end{equation}
Here $\displaystyle \boldsymbol{v}=\sum\limits_{i=1}^N\alpha_i \boldsymbol{u}_{i}$ is the average velocity of the multi-fluid, and the coefficients $\alpha_i\left(\{\xi_{j}\}_{j=1}^{N}\right)>0$ are such that $\displaystyle \sum\limits_{i=1}^N \alpha_i=1$ (for instance, $\alpha_i=\xi_i$). Thus, the values $\alpha_i$, due to their dependence on the concentrations, satisfy the equations
$$\frac{\partial \alpha_i}{\partial t}+\boldsymbol{v}\cdot\nabla\alpha_i=0,\quad i=1,\ldots,N.$$
In the simplest version, $\alpha_i$ may be regarded as constants. The main existence results for the model \eqref{mamprok1d.continit}, \eqref{mamprok1d.mominit} are obtained in \cite{mamprok1d.semi1}, \cite{mamprok1d.semi2}, \cite{mamprok1d.smz1}, \cite{mamprok1d.smz2} and \cite{mamprok1d.prokkraj}.

To conclude the Section, let us note that the multidimensional existence theorems concern only weak solutions, whose regularity is not sufficient even for the uniqueness; smoothness increase is hindered by serious obstacles. Moreover, the difficulty of multidimensional problems eclipses the study of many qualitative properties of solutions, as well as related problems including modeling; this gives researchers the right to consider the corresponding questions first in the one-dimensional case. Hence, whereas in the solvability theory for the main boundary value problems of the viscous gas there was a shift of emphasis from the one-dimensional case to the multidimensional one already two decades ago, in many other domains of the theory, the one-dimensional problems stay at the forefront.

On the other hand, as well as in the multidimensional case, the classical one-dimensional results for mono-fluids cannot be reproduced for multi-fluids automatically, in particular, due to essentially different structure of the viscous terms, namely, the presence of non-diagonal viscosity matrices; this difference in difficulty {\itshape does not depend on the
dimension of the flow}.

\section{Statement of the 1D problem, formulation of the result, Lagrangian coordinates}
Consider the rectangular $Q_{T}$ (here and below $Q_{t}=(0, 1)\times(0, t)$) with an arbitrary finite height $T>0$, and the system of equations ($i=1,\ldots,N$)
\begin{align}\label{mamprok1d.newcontinuity}\partial_{t}\rho_{i}+\partial_{x}(\rho_{i} v)=0,\quad v=\frac{1}{N}\sum\limits_{i=1}^{N}u_{i},\end{align}
\begin{align}\label{mamprok1d.newmomentum}\rho_{i}\left(\partial_{t}u_{i}+v\partial_{x}u_{i}\right)+K\partial_{x} \rho^{\gamma}=
\sum\limits_{j=1}^N \mu_{ij}\partial_{xx}u_{j},\quad \rho=\sum\limits_{i=1}^{N}\rho_{i},\end{align}
with the following initial and boundary conditions ($i=1,\ldots,N$):
\begin{align}\label{mamprok1d.nachusl}\rho_{i}|_{t=0}=\rho_{0i}(x), \quad u_{i}|_{t=0}=u_{0i}(x),\quad x\in [0, 1],\end{align}
\begin{align}\label{mamprok1d.boundvelocity}u_{i}|_{x=0}=u_{i}|_{x=1}=0,\quad t\in [0, T].\end{align}

{\bfseries Definition 1.} Refer as a {\it strong solution} to the problem \eqref{mamprok1d.newcontinuity}--\eqref{mamprok1d.boundvelocity} a collection of $2N$ functions $(\rho_{1},\ldots, \rho_{N}, u_{1},\ldots, u_{N})$ such that the equations \eqref{mamprok1d.newcontinuity},~\eqref{mamprok1d.newmomentum} are valid a.~e. in $Q_{T}$, the initial conditions \eqref{mamprok1d.nachusl}~are satisfied for a.~a. $x\in (0, 1)$, the boundary conditions \eqref{mamprok1d.boundvelocity} hold for a.~a.~$ t\in (0, T)$, and the following inequalities and inclusions hold $($$i=1,\ldots,N$$)$
$$\rho_{i}>0,\quad \rho_{i}\in L_{\infty}\big(0, T; W^{1}_{2}(0, 1)\big), \quad  \partial_{t}\rho_{i}\in L_{\infty}\big(0, T; L_{2}(0, 1)\big),$$
$$u_{i}\in L_{\infty}\big(0, T; W^{1}_{2}(0, 1)\big)\bigcap L_{2}\big(0, T; W^{2}_{2}(0, 1)\big),\quad \partial_{t}u_{i} \in L_{2}(Q_{T}).$$

The result is formulated as the following Theorem.

{\bfseries Theorem 2.} {\it Let the initial conditions in \eqref{mamprok1d.nachusl} satisfy the assumptions
$$\rho_{0i}\in W^{1}_{2}(0,1),\quad \rho_{0i}>0,\quad u_{0i}\in \overset{\circ}{W^1_2}(0, 1),\quad i=1,\ldots,N,$$
the symmetric viscosity matrix $\textbf{M}$ be positive, the adiabatic index $\gamma>1$, the constants $K,T>0$. Then there exists the unique strong solution to the problem \eqref{mamprok1d.newcontinuity}--\eqref{mamprok1d.boundvelocity} in the sense of Definition~1.}

In the model case (when all densities are equal to each other), the proof of Theorem 2 is given in \cite{mamprok1d.semi3} and \cite{mamprok1d.prok17}. Let us comment the scheme of the proof of Theorem~2. First of all, we prove the local in time solvability of the initial boundary value problem which is obtained from \eqref{mamprok1d.newcontinuity}--\eqref{mamprok1d.boundvelocity} via the Galerkin method (with respect to the spacial variable $x$) in the momentum equations~\eqref{mamprok1d.newmomentum}. The proof is made via the Leray--Schauder Fixed Point Theorem. Then, basing on the uniform estimates, we pass to the limit and establish the local solvability of the problem \eqref{mamprok1d.newcontinuity}--\eqref{mamprok1d.boundvelocity}, i.~e. on a small time interval $[0, t_{0}]$. In order to continue the solution over the whole interval $[0,T]$, we prove estimates in which the constants depend on the input data of the problem and on $T$, but not on $t_{0}$. In conclusion, we prove the uniqueness in a standard way. The key problem in the proof of Theorem 2 is the verification of the strict positiveness and boundedness of the densities (note that this problem is the key one for the complete model \eqref{mamprok1d.21061710}, \eqref{mamprok1d.21061711} as well). We consider this problem in more detail in Section 4.

During the study of the problem \eqref{mamprok1d.newcontinuity}--\eqref{mamprok1d.boundvelocity}, the parallel use of the Lagrangian coordinates is convenient. Let us accept $\displaystyle y(x,t)=\int\limits_{0}^{x}\rho(s,t)\,ds$ and~$t$ as new independent variables. Then the system \eqref{mamprok1d.newcontinuity}, \eqref{mamprok1d.newmomentum} takes the form
$$\partial_{t}\rho_{i}+\rho\rho_{i}\partial_{y}v=0,\quad v=\frac{1}{N}\sum\limits_{i=1}^{N}u_{i},$$
$$\frac{\rho_{i}}{\rho}\partial_{t}u_{i}+K\partial_{y} \rho^{\gamma}=\sum\limits_{j=1}^N \mu_{ij}\partial_{y}(\rho\partial_{y}u_{j}),\quad i=1,\ldots,N,\quad \rho=\sum\limits_{i=1}^{N}\rho_{i}.$$

The domain $Q_{T}$ is transformed into the rectangular $\Pi_{T}=(0, d)\times(0, T)$, where $\displaystyle d=\int\limits_{0}^{1}\rho_{0}\,dx>0$, $\displaystyle \rho_{0}=\sum\limits_{i=1}^{N}\rho_{0i}$, and the initial and boundary conditions take the form ($i=1,\ldots,N$)
$$\rho_{i}|_{t=0}=\widetilde{\rho}_{0i}(y), \quad u_{i}|_{t=0}=\widetilde{u}_{0i}(y),\quad y\in [0, d],$$
$$u_{i}|_{y=0}=u_{i}|_{y=d}=0,\quad t\in [0, T].$$

\section{Review of 1D viscous gas theory, strict positiveness and boundedness of the densities in the complete model \eqref{mamprok1d.21061710}, \eqref{mamprok1d.21061711}}

Consider the initial boundary value problem for the 1D Navier--Stokes equations (in the Lagrangian coordinates):
\begin{equation}\label{mamprok1d.1newcontinuity1lagr1d}\partial_{t}\rho+\rho^{2}\partial_{y}v=0,\end{equation}
\begin{equation}\label{mamprok1d.1newmomentum1lagr1d}\partial_{t}v+K\partial_{y} \rho^{\gamma}=\mu\partial_{y}(\rho\partial_{y}v),\end{equation}
\begin{equation}\label{mamprok1d.1nachusl1lagr1d}\rho|_{t=0}=\widetilde{\rho}_{0}, \quad v|_{t=0}=\widetilde{v}_{0},\quad y\in [0, d],\end{equation}
\begin{equation}\label{mamprok1d.1boundvelocity1lagr1d}v|_{y=0}=v|_{y=d}=0,\quad t\in [0, T].\end{equation}
The first a priori estimate is obtained via the standard arguments and takes the form:
\begin{equation}\label{mamprok1d.lemma1lagr1d}
\|v\|_{L_{\infty}\big(0, T;L_{2}(0, d)\big)}+\|\sqrt{\rho}\partial_{y}v\|_{L_{2}(\Pi_{T})}+\|\rho\|_{L_{\infty}\big(0, T;L_{\gamma-1}(0,d)\big)}\leqslant C.
\end{equation}
The equation \eqref{mamprok1d.1newcontinuity1lagr1d} allows to express $\rho\partial_{y}v=-\partial_{t}\ln\rho$ and to substitute it into~\eqref{mamprok1d.1newmomentum1lagr1d}:
\begin{equation}\label{mamprok1d.eq01021711d}\partial_{ty}\ln\rho+K\partial_{y}\rho^{\gamma}=-\partial_{t}v.\end{equation}
Let us multiply this equality by $\displaystyle \partial_{y}\ln\rho=:w$ and integrate over $y\in(0,d)$, then we obtain
\begin{equation}\label{mamprok1d.eq01021721d}
\frac{1}{2}\frac{d}{dt}\left(\int\limits_{0}^{d}w^{2}\, dy\right)+K\gamma\int\limits_{0}^{d}\rho^{\gamma}w^{2}\,dy=-\int\limits_{0}^{d}\left(\partial_{t} v\right)w\,dy.
\end{equation}
The right-hand side may be transformed via integration by parts due to \eqref{mamprok1d.eq01021711d}:
\begin{equation}\label{mamprok1d.eq010217771d}-\int\limits_{0}^{d}\left(\partial_{t} v\right)w\,dy=-\frac{d}{dt}\left(\int\limits_{0}^{d} vw\,dy\right)+\int\limits_{0}^{d} \rho|\partial_{y}v|^{2}\,dy.\end{equation}
Thus, after integration of \eqref{mamprok1d.eq01021721d} in $t$, using \eqref{mamprok1d.eq010217771d}, we find that
$$\|w\|^{2}_{L_{2}(0, d)}+2K\gamma\int\limits_{0}^{t}\int\limits_{0}^{d}\rho^{\gamma}w^{2}\,dyd\tau \leqslant$$
$$\leqslant\|w_{0}\|^{2}_{L_{2}(0, d)}-2\int\limits_{0}^{d}vw\,dy+2\int\limits_{0}^{d}v_{0}w_{0}\,dy+2\int\limits_{0}^{t}\|\sqrt{\rho}\partial_{y}v\|_{L_{2}(0, d)}\,d\tau,$$
where $w_{0}=w(0,t)$. Using Cauchy's inequality and the estimate \eqref{mamprok1d.lemma1lagr1d}, we derive
$$\|w(t)\|_{L_{2}(0,d)}\leqslant C\quad \forall\, t\in[0,T],$$
i.~e. the norm of the derivative $\partial_{y}\ln\rho$ in $L_{2}(0,d)$ is bounded uniformly in $t\in[0,T]$.

The equation \eqref{mamprok1d.1newcontinuity1lagr1d} and the conditions \eqref{mamprok1d.1nachusl1lagr1d}, \eqref{mamprok1d.1boundvelocity1lagr1d} obviously imply that for every $t\in[0,T]$ the equality $\rho(z(t), t)=d$ holds with some $z(t)\in[0, d]$. Hence, we can use the presentation
$${\ln{\rho(y,t)}}={\ln{\rho(z(t),t)}}+\int\limits_{z(t)}^{y}\partial_{s}\ln{\rho(s,t)}\, ds.$$
This leads to the fact that $|\ln\rho(y,t)|\leqslant |\ln{d}|+\sqrt{d}\|w\|_{L_{2}(0,d)}\leqslant C$, and hence
$$0<C^{-1}\leqslant\rho(y,t)\leqslant C.$$

For the complete multi-fluid model \eqref{mamprok1d.21061710}, \eqref{mamprok1d.21061711}, difficulties appear even on the stage of introducing the Lagrangian coordinates. In fact, if the Lagrangian coordinates are related to the velocity of one constituent (as it was made in \cite{mamprok1d.kazhpetr78}), say, with the number~$m$, then we come to the following system of equations $(i=1,\ldots,N)$:
$$\partial_{t}\rho_{i}+\rho_{m}(u_{i}-u_{m})\partial_{y}\rho_{i}+\rho_{i}\rho_{m}\partial_{y}u_{i}=0,$$
$$\frac{\rho_{i}}{\rho_{m}}\partial_{t}u_{i}+\rho_{i}(u_{i}-u_{m})\partial_{y}u_{i}+K\partial_{y} \rho^{\gamma}=\sum\limits_{j=1}^N \mu_{ij}\partial_{y}(\rho_{m}\partial_{y}u_{j})+
  \frac{1}{\rho_{m}}\sum\limits_{j=1}^{N}a_{ij}(u_{j}-u_{i}).$$
The analysis of this system via 1D viscous gas technique is possible only in the case of the diagonal viscosity matrix.

For the modified multi-fluid model \eqref{mamprok1d.newcontinuity}, \eqref{mamprok1d.newmomentum}, the direct application of the technique also fails. However, as we have seen, it is possible to overcome all difficulties and to obtain all necessary estimates.

\section{Open problems}

Let us identify several open problems in the 1D theory of viscous compressible multi-fluids:
\begin{itemize}
  \item General model with complete viscosity matrix.
  \item Heat-conductive models.
  \item Asymptotics as $t\to+\infty$.
\item Other boundary value problems.
\end{itemize}

\end{document}